\documentclass{amsart}%
\usepackage{amsmath}
\usepackage{graphicx}%
\usepackage{amsfonts}%
\usepackage{amssymb}
\newtheorem{theorem}{Theorem}
\theoremstyle{plain}

\newtheorem{corollary}[theorem]{Corollary}

\newtheorem{definition}[theorem]{Definition}

\newtheorem{lemma}[theorem]{Lemma}

\newtheorem{proposition}[theorem]{Proposition}
\newtheorem{remark}[theorem]{Remark}

\begin{document}
\title{The $\ell^{1}$-index of Tsirelson type spaces}
\author{Denny H. Leung}
\address{Department of Mathematics \\
National University of Singapore \\
Singapore 117543}
\email{matlhh@nus.edu.sg}
\author{Wee-Kee Tang}
\address{Mathematics and Mathematics Education\\
National Institute of Education\\
Nanyang Technological University\\
1 Nanyang Walk, Singapore 637616}
\email{wktang@nie.edu.sg}

\begin{abstract}
If $\alpha$ and $\beta$ are countable ordinals such that $\beta\neq0$, denote
by $\overset{_{\sim}}{T}_{\alpha,\beta}$ the completion of $c_{00}$ with
respect to the implicitly defined norm
\[
\Vert x\Vert=\max\{\Vert x\Vert_{c_{0}},\frac{1}{2}\sup\sum_{i=1}^{j}\Vert
E_{i}x\Vert\},
\]
where the supremum is taken over all finite subsets $E_{1},\dots,E_{j}$ of
$\mathbb{N}$ such that $E_{1}<\dots<E_{j}$ and $\{\min E_{1},\dots,\min
E_{j}\}\in\mbox{$\mathcal{S}$}_{\beta}$. It is shown that the Bourgain
$\ell^{1}$-index of $\overset{_{\sim}}{T}_{\alpha,\beta}\ $is $\omega
^{\alpha+\beta\cdot\omega}$. In particular, if $\omega_{1}>\alpha
=\omega^{\alpha_{1}}\cdot m_{1}+\dots+\omega^{\alpha_{n}}\cdot m_{n}$ in
Cantor normal form and $\alpha_{n}$ is not a limit ordinal, then there exists
a Banach space whose $\ell^{1}$-index is $\omega^{\alpha}$.

\end{abstract}
\maketitle

Let $E$ be a separable Banach space not containing a copy of $\ell^{1}$. The
complexity of the $\ell^{1}(n)$'s inside $E$ may be measured by Bourgain's
$\ell^{1}$-index \cite{B} or by locating so called $\ell_{\alpha}^{1}%
$-spreading models \cite{K}. It is easy to see that the existence of
$\ell_{\alpha}^{1}$-spreading models implies a large $\ell^{1}$-index. In
general, the implication is not reversible \cite[Remark 6.6(i)]{JO}. However,
suppose that $T$ is the standard Tsirelson space constructed by Figiel and
Johnson \cite{FJ} (the dual of the original Tsirelson space \cite{T}). It is
known that there is a constant $K$ such that every normalized block basic
sequence in $T$ is $K$-equivalent to a subsequence of the unit vector basis of
$T$ (see e.g., \cite{CJT}). Using this observation, one can show that the
existence of $\ell^{1}$-block trees in $T$ with large indices leads to the
existence of large $\ell_{\alpha}^{1}$-spreading models. The result can be
used to calculate the $\ell^{1}$-index of $T$. In this paper, we show that a
similar method can be applied to certain general Tsirelson type spaces. In
particular, it is shown that if $\omega_{1}>\alpha=\omega^{\alpha_{1}}\cdot
m_{1}+\dots+\omega^{\alpha_{n}}\cdot m_{n}$ in Cantor normal form and
$\alpha_{n}$ is not a limit ordinal, then there exists a Banach space whose
$\ell^{1}$-index is $\omega^{\alpha}$. This gives a partial answer to Question
1 in \cite{JO}.

If $M$ is an infinite subset of $\mathbb{N}$, denote the set of all finite,
respectively infinite subsets of $M$ by $[M]^{<\infty},$ respectively $\left[
M\right]  $. A subset $\mbox{$\mathcal{F}$}$ of $[\mathbb{N}]^{<\infty}$ is
\emph{hereditary} if $G\in\mbox{$\mathcal{F}$}$ whenever $G\subseteq F\in
\mbox{$\mathcal{F}$}$. $\mbox{$\mathcal{F}$}$ is \emph{spreading} if whenever
$F=\{n_{1},\dots,n_{k}\}\in\mbox{$\mathcal{F}$}$ with $n_{1}<\dots<n_{k}$ and
$m_{1}<\dots<m_{k}$ satisfies $m_{i}\geq n_{i}$ for $1\leq i\leq k$ then
$\{m_{1},\dots,m_{k}\}\in\mbox{$\mathcal{F}$}$. $\mbox{$\mathcal{F}$}$ is
\emph{compact} if it is compact in the product topology in $2^{\mathbb{N}}$. A
set $\mbox{$\mathcal{F}$}$ of finite subsets of $\mathbb{N}$ is called
\emph{regular} if it has all three properties. If $E$ and $F$ are finite
subsets of $\mathbb{N}$, we write $E<F$, respectively $E\leq F$, to mean $\max
E<\min F$, respectively $\max E\leq\min F$ ($\max\emptyset=0$ and
$\min\emptyset=\infty$). We abbreviate $\{n\}<E$ and $\{n\}\leq E$ to $n<E$
and $n\leq E$ respectively. Given $\mbox{$\mathcal{F}$}\subseteq
\lbrack\mathbb{N}]^{<\infty}$, a sequence of finite subsets $\{E_{1}%
,\dots,E_{n}\}$ of $\mathbb{N}$ is said to be $\mbox{$\mathcal{F}$}%
$-admissible if $E_{1}<\dots<E_{n}$ and $\{\min E_{1},\dots,\min E_{n}%
\}\in\mbox{$\mathcal{F}$}$. If $\mbox{$\mathcal{M}$}$ and $\mbox{$\mathcal{N}%
$}$ are regular subsets of $[\mathbb{N}]^{<\infty}$, we let
\[
\mbox{$\mathcal{M}$}[\mbox{$\mathcal{N}$}]=\{\cup_{i=1}^{k}F_{i}:F_{i}\in
\mbox{$\mathcal{N}$}\text{ for all $i$ and }\{F_{1},\dots,F_{k}\}\text{$\mbox
{$\mathcal{M}$}$-admissible}\}
\]
and%

\[
\left(  \mathcal{M},\mathcal{N}\right)  =\left\{  M\cup N:M<N,M\in
\mathcal{M}\text{ and }N\in\mathcal{N}\right\}  .
\]
We also write $\mathcal{M}^{2}$ for $\left(  \mathcal{M},\mathcal{M}\right)
.$ Of primary importance are the Schreier classes as defined in \cite{AA}. Let
$\mbox{$\mathcal{S}$}_{0}=\{\{n\}:n\in\mathbb{N}\}\cup\{\emptyset\}$ and
$\mbox{$\mathcal{S}$}_{1}=\{F\subseteq\mathbb{N}:|F|\leq\min F\}$. Here $|F|$
denotes the cardinality of $F$. The higher Schreier classes are defined
inductively as follows. $\mbox{$\mathcal{S}$}_{\alpha+1}=\mbox{$\mathcal{S}$%
}_{1}[\mbox{$\mathcal{S}$}_{\alpha}]$ for all $\alpha<\omega_{1}$. If $\alpha$
is a countable limit ordinal, choose a sequence $(\alpha_{n})$ strictly
increasing to $\alpha$ and set
\[
\mbox{$\mathcal{S}$}_{\alpha}=\{F:F\in\mbox{$\mathcal{S}$}_{\alpha_{n}}\text{
for some $n\leq|F|$}\}.
\]
It is clear that $\mbox{$\mathcal{S}$}_{\alpha}$ is a regular family for all
$\alpha<\omega_{1}$. If $M=(m_{1},m_{2},\dots)$ is a subsequence of
$\mathbb{N}$, let $\mbox{$\mathcal{S}$}_{\alpha}(M)=\{\{m_{i}:i\in
F\}:F\in\mbox{$\mathcal{S}$}_{\alpha}\}$. Since $\mbox{$\mathcal{S}$}_{\alpha
}$ is spreading, $\mbox{$\mathcal{S}$}_{\alpha}(M)\subseteq\mbox{$\mathcal{S}%
$}_{\alpha}$.

Let $c_{00}$ be the space of all finitely supported sequences. If $F\in\left[
\mathbb{N}\right]  ^{<\infty}$ and $a=\left(  a_{n}\right)  \in c_{00}$, let
$Fa=\left(  b_{n}\right)  \in c_{00},$ where $b_{n}=a_{n}$ if $n\in F$ and $0$
otherwise; also set $\sigma_{F}((a_{n}))=\sum_{n\in F}|a_{n}|$. Finally, if
$\mbox{$\mathcal{S}$}_{0}\subseteq\mbox{$\mathcal{F}$}\subseteq\lbrack
\mathbb{N}]^{<\infty}$, define an associated norm $\Vert\cdot\Vert
_{{\tiny \mbox{$\mathcal{F}$}}}$ on $c_{00}$ by $\Vert(a_{n})\Vert
_{{\tiny \mbox{$\mathcal{F}$}}}=\sup_{F\in{\tiny \mbox{$\mathcal{F}$}}}%
\sigma_{F}((a_{n}))$.

\begin{definition}
Let $\alpha,\beta$ be countable ordinals such that $\beta\neq0.$ Define
$\left\|  \cdot\right\|  _{n}$ and $\left\|  \cdot\right\|  _{n}^{\prime},$
$n\in\mathbb{N}$, inductively as follows. Let $\left\|  \cdot\right\|
_{0}=\left\|  \cdot\right\|  _{0}^{\prime}=\left\|  \cdot\right\|
_{\mathcal{S}_{\alpha}}.$ If $x\in c_{00},$ set
\[
\left\|  x\right\|  _{n+1}=\max\left\{  \left\|  x\right\|  _{n}%
,\,\sup\left\{  \dfrac{1}{2}\sum_{i=1}^{j}\left\|  E_{i}x\right\|
_{n}:\left\{  E_{1},...,E_{j}\right\}  \,\,\mathcal{S}_{\beta}%
\text{-admissible}\right\}  \right\}  ,
\]
and%
\[
\left\|  x\right\|  _{n+1}^{\prime}=\max\left\{  \left\|  x\right\|
_{n}^{\prime},\,\sup\left\{  \dfrac{1}{2}\sum_{i=1}^{j}\left\|  E_{i}%
x\right\|  _{n}^{\prime}:\left\{  E_{1},...,E_{j}\right\}  \,\,\left(
\mathcal{S}_{\beta}\right)  ^{2}\text{-admissible}\right\}  \right\}  .
\]
\end{definition}

Note that $\left(  \left\|  x\right\|  _{n}\right)  _{n\in\mathbb{N}}$ and
$\left(  \left\|  x\right\|  _{n}^{\prime}\right)  _{n\in\mathbb{N}}$ are
increasing sequences majorized by the $\ell^{1}$-norm of $x.$ Let $\left\|
x\right\|  _{\overset{\thicksim}{T}}=\lim\limits_{n\rightarrow\infty}\left\|
x\right\|  _{n}$ and $\left\|  x\right\|  _{\overset{\thickapprox}{T}}%
=\lim\limits_{n\rightarrow\infty}\left\|  x\right\|  _{n}^{\prime}.$ Denote by
$\overset{_{\sim}}{T}_{\alpha,\beta}$ and $\overset{\approx}{T}_{\alpha,\beta
}$ respectively the completion of $c_{00}$ under the norms $\left\|
\cdot\right\|  _{\overset{_{\sim}}{T}}$ and $\left\|  \cdot\right\|
_{\overset{\approx}{T}}.$ Clearly, $\overset{_{\sim}}{T}_{0,1}$ is the
Tsirelson space constructed by Figiel and Johnson \cite{FJ} and $\overset
{_{\sim}}{T}_{0,\beta}$ is the space denoted by $T\left(  \mathcal{S}_{\beta
},\frac{1}{2}\right)  $ in \cite{JO}. The $\ell^{1}$-index of $\overset
{_{\sim}}{T}_{0,\beta}$ is shown to be $\omega^{\beta\cdot\omega}$ in
\cite{JO}. Here, we use a different argument to compute the $\ell^{1}$-indices
of the spaces $\overset{_{\sim}}{T}_{\alpha,\beta}.$ The next proposition can
be verified immediately.

\begin{proposition}
\label{L1}The norms $\left\|  \cdot\right\|  _{\overset{_{\sim}}{T}}$ and
$\left\|  \cdot\right\|  _{\overset{\approx}{T}}$ satisfy the implicit
equations%
\[
\left\|  x\right\|  _{\overset{_{\sim}}{T}}=\max\left\{  \left\|  x\right\|
_{\mathcal{S}_{\alpha}},\,\sup\left\{  \dfrac{1}{2}\sum_{i=1}^{j}\left\|
E_{i}x\right\|  _{\overset{_{\sim}}{T}}:\left\{  E_{1},...,E_{j}\right\}
\,\mathcal{S}_{\beta}\text{-admissible}\right\}  \right\}
\]
and%
\[
\left\|  x\right\|  _{\overset{\approx}{T}}=\max\left\{  \left\|  x\right\|
_{\mathcal{S}_{\alpha}},\,\sup\left\{  \dfrac{1}{2}\sum_{i=1}^{j}\left\|
E_{i}x\right\|  _{\overset{\approx}{T}}:\left\{  E_{1},...,E_{j}\right\}
\,\left(  \mathcal{S}_{\beta}\right)  ^{2}\text{-admissible}\right\}
\right\}
\]
for all $x\in c_{00}.$
\end{proposition}

Proposition \ref{L2} is a close relative of Lemma 5 in \cite{CJT}. It is the
key observation that allows us to reduce $\ell^{1}$-block trees on
$\overset{_{\sim}}{T}_{\alpha,\beta}$ to subsequences of the unit vector basis
$\left(  e_{n}\right)  $ of $\overset{\approx}{T}_{\alpha,\beta}.$ The
following lemma is easily established by induction.

\begin{lemma}
\label{L1a}Suppose that $n_{1}\leq I_{1}<n_{2}\leq I_{2}<...<n_{k}\leq I_{k}$
and $\left|  I_{j}\right|  \leq2$ for $1\leq j\leq k.$ If $\left\{
n_{1},n_{2},...,n_{k}\right\}  \in\mathcal{S}_{\beta}$ for some $\beta
<\omega_{1},$ then $\bigcup_{j=1}^{k}I_{j}\in\left(  \mathcal{S}_{\beta
}\right)  ^{2}.$
\end{lemma}

Obviously, the coordinate unit vectors $\left(  u_{n}\right)  $ is a
normalized $1$-unconditional basis of $\overset{_{\sim}}{T}_{\alpha,\beta}.$
The \emph{support} of an element $x=\sum a_{n}u_{n}\in\overset{_{\sim}}%
{T}_{\alpha,\beta}$ is the set supp $x=\left\{  n:a_{n}\neq0\right\}  .$

\begin{proposition}
\label{L2}For every $\left\|  \cdot\right\|  _{\overset{_{\sim}}{T}}%
$-normalized block basis $\left(  x_{1},x_{2},...,x_{p}\right)  $ in
$\overset{_{\sim}}{T}_{\alpha,\beta},$ and all $\left(  a_{k}\right)  \in
c_{00},$
\[
\left\|  \sum_{k=1}^{p}a_{k}x_{k}\right\|  _{\overset{_{\sim}}{T}}%
\leq2\left\|  \sum_{k=1}^{p}a_{k}e_{i_{k}}\right\|  _{\overset{\approx}{T}}%
\]
for all $\left(  a_{k}\right)  \in c_{00},$ where $i_{k}=\max\,$%
\emph{supp}$\,x_{k},$ and $\left(  e_{n}\right)  $ is the unit vector basis of
$\overset{\approx}{T}_{\alpha,\beta}.$
\end{proposition}

\begin{proof}
With the notation as above, we prove by induction that $\left\|  \sum
_{k=1}^{p}a_{k}x_{k}\right\|  _{n}\leq2\left\|  \sum_{k=1}^{p}a_{k}e_{i_{k}%
}\right\|  _{n}^{\prime}$ for all $n\in\mathbb{N\cup}\left\{  0\right\}  ,$
$\left(  a_{k}\right)  \in c_{00}.$

To establish the inequality for the case $n=0,$ let $I\in\mathcal{S}_{\alpha
}.$ Define $J=\left\{  k:I\cap\text{supp\thinspace}x_{k}\neq\emptyset\right\}
.$ Then
\begin{align*}
\sigma_{I}\left(  \sum_{k=1}^{p}a_{k}x_{k}\right)   &  =\sum_{k=1}^{p}\left|
a_{k}\right|  \sigma_{I}\left(  x_{k}\right) \\
&  \leq\sum_{k\in J}\left|  a_{k}\right|  \left\|  x_{k}\right\|  _{0}\\
&  \leq\sum_{k\in J}\left|  a_{k}\right|  =\sigma_{L}\left(  \sum_{k=1}%
^{p}a_{k}e_{i_{k}}\right)  ,\,\text{where }L=\left\{  i_{k}:k\in J\right\}
,\\
&  \leq\left\|  \sum_{k=1}^{p}a_{k}e_{i_{k}}\right\|  _{0}^{\prime},\text{
since }L\in\mathcal{S}_{\alpha}.
\end{align*}

Suppose the proposition holds for some $n.$ Let $\left\{  E_{1},...,E_{q}%
\right\}  $ be $\mathcal{S}_{\beta}$-admissible. Without loss of generality,
we may assume that $E_{1},...,E_{q}$ are successive integer intervals, that
for all $j,$ $E_{j}\cap$\thinspace supp$\,x_{k}\neq\emptyset$ for some $k,$
and that $i_{p}\leq\max E_{q}.$ Also let $I_{k}$ be the integer interval
$\left[  i_{k-1}+1,i_{k}\right]  \,$($i_{0}\equiv0$). Let $A=\left\{
j:E_{j}\nsubseteq I_{k}\text{ for any }k\right\}  $ and $B=\left\{  j:j\notin
A\right\}  .$ For $j\in A,$ set $H_{j}=\left\{  k:I_{k}\subseteq
E_{j}\right\}  $ and $G_{j}=\left\{  i_{k}:k\in H_{j}\right\}  .$ Then define
$F_{j}=\left(  E_{j}\cap\left\{  i_{1},...,i_{p}\right\}  \right)  \setminus
G_{j}.$ Note that $F_{j}<G_{j}$ for all $j\in A.$ If $j\in B,$ set
$G_{j}=E_{j}\cap\left\{  i_{1},...,i_{p}\right\}  .$

It follows from Lemma \ref{L1a} that $\left(  F_{j}\right)  _{j\in A}%
\cup\left(  G_{j}\right)  _{j=1}^{q}$ is $\left(  \mathcal{S}_{\beta}\right)
^{2}$-admissible. Finally, $\ $let $J=\left\{  k:k\notin\cup_{j\in A}%
H_{j},\text{ }I_{k}\cap\left(  \bigcup E_{j}\right)  \neq\emptyset\right\}  .$
Now
\begin{align*}
\sum_{j=1}^{q}\left\|  E_{j}\left(  \sum_{k=1}^{p}a_{k}x_{k}\right)  \right\|
_{n}  &  =\sum_{j=1}^{q}\left\|  E_{j}\left(  \sum_{j^{\prime}\in A}\sum_{k\in
H_{j^{\prime}}}a_{k}x_{k}+\sum_{k\in J}a_{k}x_{k}\right)  \right\|  _{n}\\
&  \leq\sum_{j=1}^{q}\left(  \left\|  E_{j}\left(  \sum_{j^{\prime}\in A}%
\sum_{k\in H_{j^{\prime}}}a_{k}x_{k}\right)  \right\|  _{n}+\left\|
E_{j}\left(  \sum_{k\in J}a_{k}x_{k}\right)  \right\|  _{n}\right) \\
&  =\sum_{j\in A}\left\|  E_{j}\left(  \sum_{k\in H_{j}}a_{k}x_{k}\right)
\right\|  _{n}+\sum_{j=1}^{q}\left\|  E_{j}\left(  \sum_{k\in J}a_{k}%
x_{k}\right)  \right\|  _{n}\\
&  \leq\sum_{j\in A}\left\|  \sum_{k\in H_{j}}a_{k}x_{k}\right\|  _{n}%
+\sum_{k\in J}\left|  a_{k}\right|  \sum_{j=1}^{q}\left\|  E_{j}x_{k}\right\|
_{n}\\
&  \leq\sum_{j\in A}\left\|  \sum_{k\in H_{j}}a_{k}x_{k}\right\|  _{n}%
+2\sum_{k\in J}\left|  a_{k}\right|  \left\|  x_{k}\right\|  _{n+1}\\
&  \leq2\left(  \sum_{j\in A}\left\|  \sum_{k\in H_{j}}a_{k}e_{i_{k}}\right\|
_{n}^{\prime}+\sum_{k\in J}\left|  a_{k}\right|  \right)  \text{ by inductive
hypothesis,}\\
&  =2\left(  \sum_{j\in A}\left\|  G_{j}\left(  \sum_{k=1}^{p}a_{k}e_{i_{k}%
}\right)  \right\|  _{n}^{\prime}+\sum_{k\in J}\left|  a_{k}\right|  \right)
.
\end{align*}
If $k\in J,$ then either $\left\{  i_{k}\right\}  =F_{j}$ for some $j\in A$ or
$\left\{  i_{k}\right\}  =G_{j}$ for some $j\in B.$ Therefore%
\[
\sum_{k\in J}\left|  a_{k}\right|  \leq\sum_{j\in A}\left\|  F_{j}\left(
\sum_{k=1}^{p}a_{k}e_{i_{k}}\right)  \right\|  _{n}^{\prime}+\sum_{j\in
B}\left\|  G_{j}\left(  \sum_{k=1}^{p}a_{k}e_{i_{k}}\right)  \right\|
_{n}^{\prime}.
\]
Hence%

\begin{align*}
\sum_{j=1}^{q}\left\|  E_{j}\left(  \sum_{k=1}^{p}a_{k}x_{k}\right)  \right\|
_{n}  &  \leq2\sum_{j\in A}\left\|  G_{j}\left(  \sum_{k=1}^{p}a_{k}e_{i_{k}%
}\right)  \right\|  _{n}^{\prime}+2\sum_{j\in A}\left\|  F_{j}\left(
\sum_{k=1}^{p}a_{k}e_{i_{k}}\right)  \right\|  _{n}^{\prime}\\
&  +2\sum_{j\in B}\left\|  G_{j}\left(  \sum_{k=1}^{p}a_{k}e_{i_{k}}\right)
\right\|  _{n}^{\prime}\\
&  =2\left(  \sum_{j\in A}\left\|  F_{j}\left(  \sum_{k=1}^{p}a_{k}e_{i_{k}%
}\right)  \right\|  _{n}^{\prime}+\sum_{j=1}^{q}\left\|  G_{j}\left(
\sum_{k=1}^{p}a_{k}e_{i_{k}}\right)  \right\|  _{n}^{\prime}\right) \\
&  \leq4\left\|  \sum_{k=1}^{p}a_{k}e_{i_{k}}\right\|  _{n+1}^{\prime},\text{
as }\left(  F_{j}\right)  _{j\in A}\cup\left(  G_{j}\right)  _{j=1}^{q}\text{
is }\left(  \mathcal{S}_{\beta}\right)  ^{2}\text{-admissible.}%
\end{align*}
Thus
\[
\frac{1}{2}\sum_{j=1}^{q}\left\|  E_{j}\left(  \sum_{k=1}^{p}a_{k}%
x_{k}\right)  \right\|  _{n}\leq2\left\|  \sum_{k=1}^{p}a_{k}e_{i_{k}%
}\right\|  _{n+1}^{\prime}%
\]
whenever $\left\{  E_{1},...,E_{q}\right\}  $ is $\mathcal{S}_{\beta}%
$-admissible. It follows that%
\[
\left\|  \sum_{k=1}^{p}a_{k}x_{k}\right\|  _{n+1}\leq2\left\|  \sum_{k=1}%
^{p}a_{k}e_{i_{k}}\right\|  _{n+1}^{\prime}.
\]
This completes the induction.
\end{proof}

Let $\alpha,\beta$ be countable ordinals. Define the families $\left(
\mathcal{F}_{n}\right)  ,$ $\left(  \mathcal{F}_{n}^{\prime}\right)  ,$
$\left(  \mathcal{G}_{n}\right)  $ and $\left(  \mathcal{G}_{n}^{\prime
}\right)  $ inductively as follows: $\mathcal{F}_{0}=\mathcal{F}_{0}^{\prime
}=\mathcal{S}_{\alpha},$ $\mathcal{G}_{1}=\mathcal{S}_{\beta},$ $\mathcal{G}%
_{1}^{\prime}=\left(  \mathcal{S}_{\beta}\right)  ^{2},$ for all
$n\in\mathbb{N},$%
\[
\mathcal{F}_{n+1}=\mathcal{S}_{\beta}\left[  \mathcal{F}_{n}\right]  ,\text{
}\mathcal{F}_{n+1}^{\prime}=\left(  \mathcal{S}_{\beta}\right)  ^{2}\left[
\mathcal{F}_{n}^{\prime}\right]  \text{, }\mathcal{G}_{n+1}=\mathcal{S}%
_{\beta}\left[  \mathcal{G}_{n}\right]  \text{ and, }\mathcal{G}_{n+1}%
^{\prime}=\left(  \mathcal{S}_{\beta}\right)  ^{2}\left[  \mathcal{G}%
_{n}^{\prime}\right]  \text{.}%
\]

It is easily verified that $\mathcal{G}_{n}\left[  \mathcal{S}_{\alpha
}\right]  =\mathcal{F}_{n},$ $\mathcal{G}_{n}^{\prime}\left[  \mathcal{S}%
_{\alpha}\right]  =\mathcal{F}_{n}^{\prime},$ $\mathcal{G}_{n}\left[
\mathcal{S}_{\beta}\right]  =\mathcal{G}_{n+1}$ and $\mathcal{G}_{n}^{\prime
}\left[  \left(  \mathcal{S}_{\beta}\right)  ^{2}\right]  =\mathcal{G}%
_{n+1}^{\prime}$ for all $n\in\mathbb{N}.$ For each $n\in\mathbb{N},$ denote
the norms $\left\|  \cdot\right\|  _{\mathcal{F}_{n}}$ and $\left\|
\cdot\right\|  _{\mathcal{F}_{n}^{\prime}}$ by $\rho_{n}$ and $\rho
_{n}^{\prime}$ respectively.

\begin{proposition}
\label{L3}For all $a\in c_{00},$ and all $n\in\mathbb{N}\cup\left\{
0\right\}  ,$ $\left\|  a\right\|  _{\overset{_{\sim}}{T}}\geq\frac{1}{2^{n}%
}\rho_{n}\left(  a\right)  .$
\end{proposition}

\begin{proof}
The proof is by induction on $n.$ The case $n=0$ is clearly true by
definition. Suppose the result holds for some $n.$ Let $E\in\mathcal{F}%
_{n+1}.$ Then $E=\bigcup_{i=1}^{j}E_{i},$ where $E_{1},...,E_{j}\in
\mathcal{F}_{n},$ $E_{1}<...<E_{j},$ and $\left\{  E_{1},...,E_{j}\right\}  $
is $\mathcal{S}_{\beta}$-admissible. For any $a=\left(  a_{k}\right)  \in
c_{00},$%
\[
\sum_{k\in E}\left|  a_{k}\right|  =\sum_{i=1}^{j}\sum_{k\in E_{i}}\left|
a_{k}\right|  \leq\sum_{i=1}^{j}\rho_{n}\left(  E_{i}a\right)  \leq2^{n}%
\sum_{i=1}^{j}\left\|  E_{i}a\right\|  _{\overset{_{\sim}}{T}}\leq
2^{n+1}\left\|  a\right\|  _{\overset{_{\sim}}{T}}.
\]
Since $E\in\mathcal{F}_{n+1}$ is arbitrary, the result follows.
\end{proof}

\begin{proposition}
\label{L4}For all $a\in c_{00},$ and all $n\in\mathbb{N}\cup\left\{
0\right\}  ,$%
\[
\left\|  a\right\|  _{\overset{\approx}{T}}\leq\sum_{i=0}^{n}\frac{\rho
_{i}^{\prime}\left(  a\right)  }{2^{i}}+\frac{1}{2^{n+1}}\sup\left\{
\sum_{i=1}^{j}\left\|  E_{i}a\right\|  _{\overset{\approx}{T}}:\left\{
E_{1},...,E_{j}\right\}  \text{ }\mathcal{G}_{n+1}^{\prime}\text{-admissible}%
\right\}  .
\]
\end{proposition}

\begin{proof}
The proof is by induction on $n.$ The case $n=0$ holds by Proposition
\ref{L2}. Assume the result holds for some $n.$ Let $a\in c_{00}.$ Suppose
$\left\{  E_{1},...,E_{j}\right\}  $ is $\mathcal{G}_{n+1}^{\prime}%
$-admissible. Let $I=\left\{  i:\left\|  E_{i}a\right\|  _{\overset{\approx
}{T}}=\rho_{0}^{\prime}\left(  E_{i}a\right)  \right\}  $ and $J=\left\{
1,2,...,j\right\}  \setminus I.$ For each $i\in I,$ choose $D_{i}\subseteq
E_{i},$ $D_{i}\in\mathcal{S}_{\alpha},$ such that $\rho_{0}^{\prime}\left(
E_{i}a\right)  =\sum_{k\in D_{i}}\left|  a_{k}\right|  .$ Now $D=\cup_{i\in
I}D_{i}\in\mathcal{G}_{n+1}^{\prime}\left[  \mathcal{S}_{\alpha}\right]
=\mathcal{F}_{n+1}^{\prime}.$ Hence
\begin{equation}
\sum_{i\in I}\left\|  E_{i}a\right\|  _{\overset{\approx}{T}}=\sum_{k\in
D}\left|  a_{k}\right|  \leq\rho_{n+1}^{\prime}\left(  a\right)  . \label{1}%
\end{equation}
On the other hand, for each $i\in J,$ there exist $\left(  \mathcal{S}_{\beta
}\right)  ^{2}$-admissible sets $\left\{  E_{1}^{i},...,E_{k_{i}}^{i}\right\}
,$ $E_{1}^{i}\cup...\cup E_{k_{i}}^{i}\subseteq E_{i}$ such that
\[
\left\|  E_{i}a\right\|  _{\overset{\approx}{T}}=\frac{1}{2}\sum_{p=1}^{k_{i}%
}\left\|  E_{p}^{i}a\right\|  _{\overset{\approx}{T}}.
\]
Now $\left\{  \min E_{p}^{i}:i\in J,\text{ }1\leq p\leq k_{i}\right\}
\in\mathcal{G}_{n+1}^{\prime}\left[  \left(  \mathcal{S}_{\beta}\right)
^{2}\right]  =\mathcal{G}_{n+2}^{\prime}.$ Hence $\left(  E_{p}^{i}\right)
_{i\in J,\text{ }1\leq p\leq k_{i}}$ is $\mathcal{G}_{n+2}^{\prime}$
admissible. Thus
\begin{align}
\sum_{i\in J}\left\|  E_{i}a\right\|  _{\overset{\approx}{T}}  &  =\frac{1}%
{2}\sum_{i\in J}\sum_{p=1}^{k_{i}}\left\|  E_{p}^{i}a\right\|  _{\overset
{\approx}{T}}\label{2}\\
&  \leq\frac{1}{2}\sup\left\{  \sum_{i=1}^{\ell}\left\|  F_{i}a\right\|
_{\overset{\approx}{T}}:\left\{  F_{1},...,F_{\ell}\right\}  \text{
}\mathcal{G}_{n+2}^{\prime}\text{-admissible}\right\}  .\nonumber
\end{align}
From the inductive hypothesis and inequalities (\ref{1}) and (\ref{2}) we get
\begin{align*}
\left\|  a\right\|  _{\overset{\approx}{T}}  &  \leq\sum_{i=0}^{n}\frac
{\rho_{i}^{\prime}\left(  a\right)  }{2^{i}}+\frac{1}{2^{n+1}}\left(
\rho_{n+1}^{\prime}\left(  a\right)  +\frac{1}{2}\sup\left\{  \sum_{i=1}%
^{\ell}\left\|  F_{i}a\right\|  _{\overset{\approx}{T}}\right\}  \right) \\
&  =\sum_{i=0}^{n+1}\frac{\rho_{i}^{\prime}\left(  a\right)  }{2^{i}}+\frac
{1}{2^{n+2}}\sup\left\{  \sum_{i=1}^{\ell}\left\|  F_{i}a\right\|
_{\overset{\approx}{T}}\right\}  ,
\end{align*}
where both suprema are taken over all $\mathcal{G}_{n+2}^{\prime}$-admissible
sets $\left\{  F_{1},...,F_{\ell}\right\}  .$ This completes the induction.
\end{proof}

Endow $\left[  \mathbb{N}\right]  ^{<\infty}$ with the product topology
inherited from $2^{\mathbb{N}}.$ If $\mathcal{F}$ is a closed subset of
$\left[  \mathbb{N}\right]  ^{<\infty},$ let $\mathcal{F}^{\prime}$ be the set
of all limit points of $\mathcal{F}.$ Define a transfinite sequence of sets
$\left(  \mathcal{F}^{\left(  \alpha\right)  }\right)  _{\alpha<\omega_{1}}$
as follows: $\mathcal{F}^{\left(  0\right)  }=\mathcal{F},$ $\mathcal{F}%
^{\left(  \alpha+1\right)  }=\left(  \mathcal{F}^{\left(  \alpha\right)
}\right)  ^{\prime}$ for all $\alpha<\omega_{1};$ $\mathcal{F}^{\left(
\alpha\right)  }=\bigcap_{\beta<\alpha}\mathcal{F}^{\left(  \beta\right)  }$
if $\alpha$ is a countable limit ordinal.

\begin{definition}
\emph{(\cite{OTW}) }Let $\mathcal{F}\subseteq\left[  \mathbb{N}\right]
^{<\infty}$ be regular$.$ Define $\iota\left(  \mathcal{F}\right)  $ to be the
unique countable ordinal $\alpha$ such that $\mathcal{F}^{\left(
\alpha\right)  }=\left\{  \emptyset\right\}  .$
\end{definition}

Let $F\in\mathcal{M}\left[  \mathcal{N}\right]  ,$ $F\neq\emptyset,$ where
$\mathcal{M}$ and $\mathcal{N}$ are regular families$.$ There exists a largest
$k\in F$ such that $F\cap\left[  1,k\right]  \in\mathcal{N}$. Set $F_{1}%
=F\cap\left[  1,k\right]  .$ If $F_{1},F_{2},...,F_{n-1}$ have been defined
and $F\setminus\bigcup_{i=1}^{n-1}F_{i}\neq\emptyset,$ set $F_{n}=\left(
F\setminus\bigcup_{i=1}^{n-1}F_{i}\right)  \cap\left[  1,k^{\prime}\right]  ,$
where $k^{\prime}$ is the largest integer in $F$ such that $\left(
F\setminus\bigcup_{i=1}^{n-1}F_{i}\right)  \cap\left[  1,k^{\prime}\right]
\in\mathcal{N}$. Since $F$ is finite, there exists an $n$ such that
$F=\bigcup_{i=1}^{n}F_{i}.$ Now $F\in\mathcal{M}\left[  \mathcal{N}\right]  $
implies that there exists an $\mathcal{M}$-admissible collection $\left\{
G_{1},...,G_{m}\right\}  $ such that $F=\bigcup_{j=1}^{m}G_{j}$ and $G_{j}%
\in\mathcal{N}$, $1\leq j\leq m.$ By the choice of $F_{i},$ and the fact that
$\mathcal{N}$ is hereditary, it is easy to see that $\min G_{i}\leq\min
F_{i},$ $1\leq i\leq n.$ Thus $\left\{  F_{1},...,F_{n}\right\}  $ is
$\mathcal{M}$-admissible, as $\mathcal{M}$ is spreading. We call $\left(
F_{i}\right)  _{i=1}^{n}$ the standard representation of $F$ (as an element of
$\mathcal{M}\left[  \mathcal{N}\right]  $).

\begin{remark}
\label{R1}\emph{Suppose that }$\left(  F_{i}\right)  _{i=1}^{n}$\emph{ and
}$\left(  G_{i}\right)  _{i=1}^{m}$\emph{ are the standard representations of
}$F$\emph{ and }$G$\emph{ respectively. If }$\ell,k\in\mathcal{N}$\emph{ are
such that }$F\cap\left[  1,\ell\right]  =G\cap\left[  1,\ell\right]  $\emph{
and }$\max F_{k}\leq\ell,$\emph{ then by construction, }$F_{i}=G_{i},$\emph{
}$1\leq i\leq k.$
\end{remark}

\begin{lemma}
\label{3.1}Let $\mathcal{M}$, $\mathcal{N}\subseteq\left[  \mathbb{N}\right]
^{<\infty}$ be regular. Suppose that $\iota\left(  \mathcal{N}\right)
=\alpha,$ then $\left(  \mathcal{M}\left[  \mathcal{N}\right]  \right)
^{\left(  \alpha\right)  }=\left(  \mathcal{M}^{\left(  1\right)  }\right)
\left[  \mathcal{N}\right]  .$
\end{lemma}

\begin{proof}
Let $F\in\left(  \mathcal{M}^{\left(  1\right)  }\right)  \left[
\mathcal{N}\right]  ,$ then $F$ can be written as $F=\cup_{i=1}^{n}F_{i},$
where $F_{1}<F_{2}<...<F_{n},$ $F_{1},...,F_{n}\in\mathcal{N}$, and $\left\{
\min F_{1},...,\min F_{n}\right\}  \in\mathcal{M}^{\left(  1\right)  }.$ In
particular, there exists $k_{0}>\max F$ such that $\left\{  \min
F_{1},...,\min F_{n},k\right\}  \in\mathcal{M}$ for all $k\geq k_{0}.$
Therefore,%
\begin{equation}
\text{for all }G\in\mathcal{N},\min G\geq k_{0},\,F\cup G\in\mathcal{M}\left[
\mathcal{N}\right]  . \label{e3}%
\end{equation}
Note that as $\mathcal{N}$ is spreading,%
\begin{equation}
\iota\left(  \left\{  G\in\mathcal{N}:\min G\geq k_{0}\right\}  \right)
=\iota\left(  \mathcal{N}\right)  =\alpha. \label{e4}%
\end{equation}
From (\ref{e3}),
\[
\left(  \left\{  F\cup G:G\in\mathcal{N},\,\min G\geq k_{0}\right\}  \right)
^{\left(  \alpha\right)  }\subseteq\left(  \mathcal{M}\left[  \mathcal{N}%
\right]  \right)  ^{\left(  \alpha\right)  }.
\]
But from (\ref{e4}), $F\in\left(  \left\{  F\cup G:G\in\mathcal{N},\,\min
G\geq k_{0}\right\}  \right)  ^{\left(  \alpha\right)  }.$ Hence $F\in\left(
\mathcal{M}\left[  \mathcal{N}\right]  \right)  ^{\left(  \alpha\right)  }.$

Conversely, we prove by induction that $\left(  \mathcal{M}\left[
\mathcal{N}\right]  \right)  ^{\left(  \gamma\right)  }\subseteq\left(
\left(  \mathcal{M}^{\left(  1\right)  }\right)  \left[  \mathcal{N}\right]
,\mathcal{N}^{\left(  \gamma\right)  }\right)  $ for all $\gamma\leq\alpha.$
The cases where $\gamma=0$ is clear.

Suppose the claim is true for some $\gamma<\alpha.$ Let $F\in\left(
\mathcal{M}\left[  \mathcal{N}\right]  \right)  ^{\left(  \gamma+1\right)  }$
with standard representation $\left(  F_{i}\right)  _{i=1}^{n}$ as an element
of $\mathcal{M}\left[  \mathcal{N}\right]  .$ Choose a sequence $\left(
G_{k}\right)  $ in $\left(  \mathcal{M}\left[  \mathcal{N}\right]  \right)
^{\left(  \gamma\right)  }\subseteq\left(  \left(  \mathcal{M}^{\left(
1\right)  }\right)  \left[  \mathcal{N}\right]  ,\mathcal{N}^{\left(
\gamma\right)  }\right)  $ that converges nontrivially to $F.$ We may assume
that $G_{k}\cap\left[  1,\min F_{n}\right]  =F\cap\left[  1,\min F_{n}\right]
$ for all $k.$ Now we may write $G_{k}=P_{k}\cup Q_{k},$ where $P_{k}<Q_{k}$,
$P_{k}\in\left(  \mathcal{M}^{\left(  1\right)  }\right)  \left[
\mathcal{N}\right]  $ and $Q_{k}\in\mathcal{N}^{\left(  \gamma\right)  }.$ Let
$P=\bigcup_{i=1}^{n-1}F_{i}$. Note that $P\in\left(  \mathcal{M}^{\left(
1\right)  }\right)  \left[  \mathcal{N}\right]  .$ We consider two cases.

\textbf{Case 1. }There exists $k$ such that $\min F_{n}\in P_{k}.$

In this case, $P\cap\left[  1,\max F_{n-1}\right]  =F\cap\left[  1,\max
F_{n-1}\right]  =G_{k}\cap\left[  1,\max F_{n-1}\right]  =P_{k}\cap\left[
1,\max F_{n-1}\right]  .$ It is clear that $\left(  F_{i}\right)  _{i=1}%
^{n-1}$ is the standard representation of $P$ as an element of $\left(
\mathcal{M}^{\left(  1\right)  }\right)  \left[  \mathcal{N}\right]  .$ By
Remark \ref{R1}, the standard representation of $P_{k}$ as an element of
$\left(  \mathcal{M}^{\left(  1\right)  }\right)  \left[  \mathcal{N}\right]
$ has the form $\left(  F_{1},...,F_{n-1},P_{k}^{n},...,P_{k}^{m}\right)  .$
In particular,
\[
\left\{  \min F_{1},...,\min F_{n-1},\min F_{n}\right\}  =\left\{  \min
F_{1},...,\min F_{n-1},\min P_{k}^{n}\right\}  \in\mathcal{M}^{\left(
1\right)  }.
\]
Thus $F=\bigcup_{i=1}^{n}F_{i}\in\left(  \mathcal{M}^{\left(  1\right)
}\right)  \left[  \mathcal{N}\right]  \subseteq\left(  \left(  \mathcal{M}%
^{\left(  1\right)  }\right)  \left[  \mathcal{N}\right]  ,\mathcal{N}%
^{\left(  \gamma+1\right)  }\right)  ,$ as required.

\textbf{Case 2. }Suppose $\min F_{n}\notin P_{k}$ for all $k\in\mathbb{N}.$

In this case, $G_{k}\cap\lbrack\min F_{n},\infty)\subseteq Q_{k}$ for all $k.$
Hence $G_{k}\cap\lbrack\min F_{n},\infty)\in\mathcal{N}^{\left(
\gamma\right)  }$ for all $k.$ Furthermore, $G_{k}\cap\lbrack\min F_{n}%
,\infty)$ converges to $F\cap\lbrack\min F_{n},\infty)=F_{n}$ nontrivially.
Thus $F_{n}\in\mathcal{N}^{\left(  \gamma+1\right)  }.$ Therefore $F=P\cup
F_{n}\in\left(  \left(  \mathcal{M}^{\left(  1\right)  }\right)  \left[
\mathcal{N}\right]  ,\mathcal{N}^{\left(  \gamma+1\right)  }\right)  ,$ as required.

Suppose $\gamma\leq\alpha$ is a limit ordinal and the result holds for all
$\eta<\gamma.$ Let $F\in\left(  \mathcal{M}\left[  \mathcal{N}\right]
\right)  ^{\left(  \gamma\right)  }$ have standard representation $\left(
F_{i}\right)  _{i=1}^{n}$ as an element of $\mathcal{M}\left[  \mathcal{N}%
\right]  .$ By the inductive hypothesis, for each $\eta<\gamma,$ $F=P_{\eta
}\cup Q_{\eta},$ where $P_{\eta}<Q_{\eta},$ $P_{\eta}\in\left(  \mathcal{M}%
^{\left(  1\right)  }\right)  \left[  \mathcal{N}\right]  $ and $Q_{\eta}%
\in\mathcal{N}^{\left(  \eta\right)  }.$ By the argument in case 1 above, if
there exists $\eta$ such that $\min F_{n}\in P_{\eta},$ then $F\in\left(
\mathcal{M}^{\left(  1\right)  }\right)  \left[  \mathcal{N}\right]
\subseteq\left(  \left(  \mathcal{M}^{\left(  1\right)  }\right)  \left[
\mathcal{N}\right]  ,\mathcal{N}^{\left(  \eta\right)  }\right)  .$ Otherwise,
$F_{n}\subseteq Q_{\eta}\in\mathcal{N}^{\left(  \eta\right)  }$ for all
$\eta<\gamma.$ Hence $F=\left(  \bigcup_{i=1}^{n-1}F_{i}\right)  \cup F_{n}%
\in\left(  \left(  \mathcal{M}^{\left(  1\right)  }\right)  \left[
\mathcal{N}\right]  ,\mathcal{N}^{\left(  \gamma\right)  }\right)  .$ This
completes the induction.
\end{proof}

\begin{proposition}
\label{P7}Let $\mathcal{M}$, $\mathcal{N}\subseteq\left[  \mathbb{N}\right]
^{<\infty}$ be regular. Suppose that $\iota\left(  \mathcal{N}\right)
=\alpha.$ Then for all $\beta<\omega_{1},$ $\left(  \mathcal{M}\left[
\mathcal{N}\right]  \right)  ^{\left(  \alpha\cdot\beta\right)  }=\left(
\mathcal{M}^{\left(  \beta\right)  }\right)  \left[  \mathcal{N}\right]  .$
\end{proposition}

\begin{proof}
The proof is by induction on $\beta.$ The case $\beta=0$ is clear. Suppose the
result is true for some $\beta.$ Then
\begin{align*}
\left(  \mathcal{M}\left[  \mathcal{N}\right]  \right)  ^{\left(  \alpha
\cdot\left(  \beta+1\right)  \right)  }  &  =\left(  \mathcal{M}\left[
\mathcal{N}\right]  \right)  ^{\left(  \alpha\cdot\beta+\alpha\right)  }\\
&  =\left(  \left(  \mathcal{M}\left[  \mathcal{N}\right]  \right)  ^{\left(
\alpha\cdot\beta\right)  }\right)  ^{\left(  \alpha\right)  }\\
&  =\left(  \left(  \mathcal{M}^{\left(  \beta\right)  }\right)  \left[
\mathcal{N}\right]  \right)  ^{\left(  \alpha\right)  }\text{ by the inductive
hypothesis,}\\
&  =\left(  \left(  \mathcal{M}^{\left(  \beta\right)  }\right)  ^{\left(
1\right)  }\left[  \mathcal{N}\right]  \right)  \text{ by Lemma \ref{3.1},}\\
&  =\left(  \mathcal{M}^{\left(  \beta+1\right)  }\right)  \left[
\mathcal{N}\right]  .
\end{align*}
Suppose the proposition is true for all $\beta<\beta_{0},$ where $\beta
_{0}<\omega_{1}$ is some limit ordinal. Clearly,
\[
\left(  \mathcal{M}^{\left(  \beta_{0}\right)  }\right)  \left[
\mathcal{N}\right]  \subseteq\bigcap_{\beta<\beta_{0}}\left(  \mathcal{M}%
^{\left(  \beta\right)  }\right)  \left[  \mathcal{N}\right]  =\bigcap
_{\beta<\beta_{0}}\left(  \mathcal{M}\left[  \mathcal{N}\right]  \right)
^{\left(  \alpha\cdot\beta\right)  }=\left(  \mathcal{M}\left[  \mathcal{N}%
\right]  \right)  ^{\left(  \alpha\cdot\beta_{0}\right)  }.
\]
On the other hand, let $F\in\bigcap_{\beta<\beta_{0}}\left(  \mathcal{M}%
^{\left(  \beta\right)  }\right)  \left[  \mathcal{N}\right]  $ have standard
representation $\left(  F_{i}\right)  _{i=1}^{n}$ as an element of
$\mathcal{M}\left[  \mathcal{N}\right]  .$ It is clear that $\left(
F_{i}\right)  _{i=1}^{n}$ is also the standard representation of $F$ as an
element of $\left(  \mathcal{M}^{\left(  \beta\right)  }\right)  \left[
\mathcal{N}\right]  $ for any $\beta<\beta_{0}.$ In particular, $\left\{  \min
F_{i}:1\leq i\leq n\right\}  \in\mathcal{M}^{\left(  \beta\right)  }$ for all
$\beta<\beta_{0}.$ Hence $\left\{  \min F_{i}:1\leq i\leq n\right\}
\in\mathcal{M}^{\left(  \beta_{0}\right)  }.$ It follows that $F\in\left(
\mathcal{M}^{\left(  \beta_{0}\right)  }\right)  \left[  \mathcal{N}\right]
.$ This completes the proof.
\end{proof}

It is well known that $\iota\left(  \mathcal{S}_{\gamma}\right)
=\omega^{\gamma}$ for all $\gamma<\omega_{1}$ (\cite[Proposition 4.10]{AA}).
The indices of $\mathcal{F}_{n}$ and $\mathcal{F}_{n}^{\prime}$ can be
computed readily with the help of Proposition \ref{P7}.

\begin{corollary}
\label{C7}$\iota\left(  \mathcal{F}_{n}\right)  =\omega^{\alpha+\beta\cdot
n},$ $\iota\left(  \mathcal{F}_{n}^{\prime}\right)  =\omega^{\alpha+\beta\cdot
n}\cdot2.$
\end{corollary}

Before proceeding further, let us recall the relevant terminology concerning
trees. A \emph{tree} on a set $X$ is a subset $T$ of $\cup_{n=1}^{\infty}%
X^{n}$ such that $(x_{1},\dots,x_{n})\in T$ whenever $n\in\mathbb{N}$ and
$(x_{1},\dots,x_{n+1})\in T$. These are the only kind of trees we will
consider. A tree $T$ is \emph{well-founded} if there is no infinite sequence
$(x_{n})$ in $X$ such that $(x_{1},\dots,x_{n})\in T$ for all $n$. Given a
well-founded tree $T$, we define the \emph{derived tree} $D(T)$ to be the set
of all $(x_{1},\dots,x_{n})\in T$ such that $(x_{1},\dots,x_{n},x)\in T$ for
some $x\in X$. Inductively, we let $D^{0}(T)=T$, $D^{\alpha+1}(T)=D(D^{\alpha
}(T))$, and $D^{\alpha}(T)=\cap_{\beta<\alpha}D^{\beta}(T)$ if $\alpha$ is a
limit ordinal. The \emph{order} of a well-founded tree $T$ is the smallest
ordinal $o(T)$ such that $D^{o(T)}(T)=\emptyset$. If $E$ is a Banach space and
$1\leq K<\infty$, an $\ell^{1}$-$K$ tree on $E$ is a tree $T$ on $S(X)=\{x\in
E:\Vert x\Vert=1\}$ such that $\Vert\sum_{i=1}^{n}a_{i}x_{i}\Vert\geq
K^{-1}\sum_{i=1}^{n}|a_{i}|$ whenever $(x_{1},\dots,x_{n})\in T$ and
$(a_{i})\subseteq\mathbb{R}$. If $E$ has a basis $(e_{i})$, a \emph{block
tree} on $E$ is a tree $T$ on $E$ so that every $(x_{1},\dots,x_{n})\in T$ is
a finite block basis of $(e_{i})$. An $\ell^{1}$-$K$-block tree on $E$ is a
block tree that is also an $\ell^{1}$-$K$ tree. The index $I(E,K)$ is defined
to be $\sup\{o(T):T\text{ is an $\ell^{1}$-$K$ tree on $E$}\}$. If $E$ has a
basis $(e_{i})$, the index $I_{b}(E,K)$ is defined similarly, with the
supremum taken over all $\ell^{1}$-$K$ block trees. The Bourgain $\ell^{1}%
$-\emph{index} of $E$ is the ordinal $I(E)=\sup\{I(E,K):1\leq K<\infty\}$. The
index $I_{b}(E)$ is defined similarly. Bourgain proved that if $E$ is a
separable Banach space not containing a copy of $\ell^{1}$, then
$I(E)<\omega_{1}$ \cite{B}. More on these and related indices can be found in
\cite{JO} and \cite{AJO}.

\begin{proposition}
\label{L8}Let $T$ be a well-founded block tree on some basis $\left(
e_{i}\right)  .$ Define
\[
\mathcal{F}\left(  T\right)  =\left\{  \left\{  \max\text{\emph{supp\thinspace
}}x_{i}:i=1,...,n\right\}  :\left(  x_{1},x_{2},...,x_{n}\right)  \in
T\right\}  .
\]
and
\[
\mathcal{G}\left(  T\right)  \mathcal{=}\left\{  G:\exists F\in\mathcal{F}%
\left(  T\right)  ,f:\mathbb{N\rightarrow N}\text{ strictly increasing, such
that }G\subseteq f\left(  F\right)  \right\}  .
\]
Then $\iota\left(  \mathcal{G}\left(  T\right)  \right)  \geq o\left(
T\right)  .$
\end{proposition}

\begin{proof}
Let $\xi=o\left(  T\right)  .$ The proof is by induction on $\xi.$ If
$o\left(  T\right)  =1,$ then $\mathcal{G}\left(  T\right)  \mathcal{\supseteq
}\left\{  \left\{  k\right\}  :k\geq n\right\}  $ for some $n\in\mathbb{N}.$
Therefore $\left(  \mathcal{G}\left(  T\right)  \right)  ^{\left(  1\right)
}\mathcal{\supseteq}\left\{  \emptyset\right\}  $ and hence $\iota\left(
\mathcal{G}\left(  T\right)  \right)  \geq1=o\left(  T\right)  .$

Suppose the proposition holds for some $\xi<\omega_{1}.$ Let $T$ be a well
founded block tree with $o\left(  T\right)  =\xi+1.$ For each $\left(
x\right)  \in T,$ let
\[
T_{x}=\bigcup_{n=1}^{\infty}\left\{  \left(  x_{1},...,x_{n}\right)  :\left(
x,x_{1},...,x_{n}\right)  \in T\right\}  .
\]
According to \cite[Proposition 4]{B}, $o\left(  T\right)  =\sup_{\left(
x\right)  \in T}\left\{  o\left(  T_{x}\right)  +1\right\}  .$ Therefore,
there exists $\left(  x_{0}\right)  \in T$ such that $o\left(  T_{x_{0}%
}\right)  =\xi.$ By the inductive hypothesis, $\iota\left(  \mathcal{G}\left(
T_{x_{0}}\right)  \right)  \geq\xi.$ Let $k_{0}=\max$ supp $x_{0}.$ Then
$G\in\mathcal{G}\left(  T_{x_{0}}\right)  $ implies $\left\{  k_{0}\right\}
\cup G\in\mathcal{G}\left(  T\right)  .$ Thus $\left\{  k_{0}\right\}
\in\left(  \mathcal{G}\left(  T\right)  \right)  ^{\left(  \xi\right)  }.$
Since $\left(  \mathcal{G}\left(  T\right)  \right)  ^{\left(  \xi\right)  }$
is spreading, $\left\{  k\right\}  \in\left(  \mathcal{G}\left(  T\right)
\right)  ^{\left(  \xi\right)  }$ for all $k\geq k_{0}.$ It follows that
$\emptyset\in\left(  \mathcal{G}\left(  T\right)  \right)  ^{\left(
\xi+1\right)  }.$ Hence $\iota\left(  \mathcal{G}\left(  T\right)  \right)
\geq\xi+1=o\left(  T\right)  .$

Suppose $o\left(  T\right)  =\xi_{0},$ where $\xi_{0}$ is a countable limit
ordinal and the proposition holds for all $\xi<\xi_{0}$. Since $o\left(
T\right)  >\xi$ for all $\xi<\xi_{0}.$ By the inductive hypothesis,
$\iota\left(  \mathcal{G}\left(  T\right)  \right)  >\xi$ for all $\xi<\xi
_{0}.$ Hence $\iota\left(  \mathcal{G}\left(  T\right)  \right)  \geq\xi
_{0}=o\left(  T\right)  .$ This completes the induction.
\end{proof}

If $\left(  x_{k}\right)  _{k=1}^{n}$ and $\left(  y_{k}\right)  _{k=1}^{n}$
are sequences in possibly different normed spaces, and $0<K<\infty,$ we write
$\left(  x_{k}\right)  _{k=1}^{n}\overset{K}{\succeq}\left(  y_{k}\right)
_{k=1}^{n}$ to mean $K\left\|  \sum_{k=1}^{n}a_{k}x_{k}\right\|  \geq\left\|
\sum_{k=1}^{n}a_{k}y_{k}\right\|  $ for all $\left(  a_{k}\right)
\subseteq\mathbb{R}.$

\begin{theorem}
\label{T14}$I\left(  \overset{_{\sim}}{T}_{\alpha,\beta}\right)  =I_{b}\left(
\overset{_{\sim}}{T}_{\alpha,\beta}\right)  =\omega^{\alpha+\beta\cdot\omega}.$
\end{theorem}

\begin{proof}
If $I_{b}\left(  \overset{_{\sim}}{T}_{\alpha,\beta}\right)  >\omega
^{\alpha+\beta\cdot\omega},$ then $I_{b}\left(  \overset{_{\sim}}{T}%
_{\alpha,\beta},K\right)  >\omega^{\alpha+\beta\cdot\omega}$ for some $K>1.$
Hence there exists an $\ell^{1}$-$K$-block tree $T$ on $\overset{_{\sim}}%
{T}_{\alpha,\beta}$ such that $o\left(  T\right)  =\xi>\omega^{\alpha
+\beta\cdot\omega}.$ Given $F\in\mathcal{F}\left(  T\right)  ,$ there exists
$\left(  x_{1},x_{2},...,x_{n}\right)  \in T$ such that $F=\{\max
\,$supp\emph{\thinspace}$x_{i}\}_{i=1}^{n}.$ According to Proposition
\ref{L2}, $\left(  e_{k}\right)  _{k\in F}\overset{2}{\succeq}\left(
x_{1},x_{2},...,x_{n}\right)  ,$ where $\left(  e_{k}\right)  _{k=1}^{\infty}$
is the unit vector basis of $\overset{\approx}{T}_{\alpha,\beta}.$ Since
$\left(  x_{1},x_{2},...,x_{n}\right)  \in T,$ $\left(  x_{1},x_{2}%
,...,x_{n}\right)  \overset{K}{\succeq}\ell^{1}\left(  \left|  F\right|
\right)  $-basis$.$ Therefore, $\left(  e_{k}\right)  _{k\in F}\overset
{2K}{\succeq}\ell^{1}\left(  \left|  F\right|  \right)  $-basis for all
$F\in\mathcal{F}\left(  T\right)  .$ Since it is clear that $\left\|  \sum
a_{k}e_{f(k)}\right\|  _{\overset{\approx}{T}}\geq\left\|  \sum a_{k}%
e_{k}\right\|  _{\overset{\approx}{T}}$ for all $\left(  a_{k}\right)  \in
c_{00}$ whenever $f:\mathbb{N}\rightarrow\mathbb{N}$ is strictly increasing,
it follows that $\left(  e_{k}\right)  _{k\in G}\overset{2K}{\succeq}\ell
^{1}\left(  \left|  G\right|  \right)  $-basis for all $G\in\mathcal{G}\left(
T\right)  .$ By Proposition \ref{L8}, $\left(  \mathcal{G}\left(  T\right)
\right)  ^{\left(  \omega^{\alpha+\beta\cdot\omega}+1\right)  }\neq\emptyset.$
Thus by \cite[Corollary 1.2]{G}, there exists $L\in\left[  \mathbb{N}\right]
$ such that $\mathcal{S}_{\alpha+\beta\cdot\omega}\cap\left[  L\right]
^{<\infty}\subseteq\mathcal{G}\left(  T\right)  .$ Hence, for all $\left(
a_{k}\right)  \in c_{00}$ and all $G\in\mathcal{S}_{\alpha+\beta\cdot\omega
}\cap\left[  L\right]  ^{<\infty},$%
\begin{equation}
\left\|  \sum_{k\in G}a_{k}e_{k}\right\|  _{\overset{\approx}{T}}\geq\frac
{1}{2K}\sum_{k\in G}\left|  a_{k}\right|  . \label{e5}%
\end{equation}
Choose $m\in\mathbb{N}$ such that $2^{m}>2K.$ According to Corollary \ref{C7},
$\iota\left(  \mathcal{F}_{i}^{\prime}\right)  =\omega^{\alpha+\beta\cdot
i}\cdot2$ for all $i=1,2,...,m.$ Applying \cite[Corollary 1.2]{G}, we obtain
$M\in\left[  L\right]  $ such that $\mathcal{F}_{i}^{\prime}\cap\left[
M\right]  ^{<\infty}\subseteq\mathcal{S}_{\alpha+\beta\cdot m+1}$ for all
$i=1,2,...,m.$ By \cite[Proposition 3.6]{OTW}, there exists $F\in
\mathcal{S}_{\alpha+\beta\cdot\omega}\left(  M\right)  \subseteq
\mathcal{S}_{\alpha+\beta\cdot\omega}\cap\left[  M\right]  ^{<\infty}$ and
$\left(  a_{j}\right)  _{j\in F}\subseteq\mathbb{R}^{+}$ such that $\sum
a_{j}=1$ and if $G\subseteq F$ with $G\in\mathcal{S}_{\alpha+\beta\cdot m+1},$
then $\sum_{j\in G}a_{j}<\frac{1}{8K}.$ Consider $x=\sum_{j\in F}a_{j}e_{j}%
\in\overset{\approx}{T}_{\alpha,\beta}.$ If $1\leq i\leq m$ and $I\in
\mathcal{F}_{i}^{\prime},$ then $I\cap F\in\mathcal{F}_{i}^{\prime}\cap\left[
M\right]  ^{<\infty}\subseteq\mathcal{S}_{\alpha+\beta\cdot m+1}.$ Hence
$\sigma_{I}\left(  x\right)  =\sigma_{I\cap F}\left(  x\right)  <\frac{1}%
{8K}.$ It follows that $\rho_{i}^{\prime}\left(  x\right)  \leq\frac{1}{8K}$
for $1\leq i\leq m.$ By Proposition \ref{L4},
\begin{align*}
\left\|  x\right\|  _{\overset{\approx}{T}}  &  \leq\sum_{i=0}^{m}\frac
{\rho_{i}^{\prime}\left(  x\right)  }{2^{i}}+\frac{1}{2^{m+1}}\sup\left\{
\sum_{i=1}^{j}\left\|  E_{i}x\right\|  _{\overset{\approx}{T}}:\left\{
E_{1},...,E_{j}\right\}  \text{ }\mathcal{G}_{m+1}^{\prime}\text{-admissible}%
\right\} \\
&  \leq\sum_{i=0}^{m}\frac{\frac{1}{8K}}{2^{i}}+\frac{1}{2^{m+1}}\left\|
x\right\|  _{\ell^{1}}<\frac{1}{2K},
\end{align*}
contrary to (\ref{e5}). This proves that $I_{b}\left(  \overset{_{\sim}}%
{T}_{\alpha,\beta}\right)  \leq\omega^{\alpha+\beta\cdot\omega}.$ On the other
hand, according to Proposition \ref{L3}, for any $n\in\mathbb{N}$, $\left\|
a\right\|  _{\overset{_{\sim}}{T}}\geq\frac{1}{2^{n}}\left\|  a\right\|
_{\mathcal{F}_{n}}$ for any $a\in c_{00}.$ By Corollary \ref{C7},
$\iota\left(  \mathcal{F}_{n}\right)  =\omega^{\alpha+\beta\cdot n}.$
Therefore, there exists an $\ell^{1}$-$2^{n}-$block basis tree $T_{n}$ on
$\overset{_{\sim}}{T}_{\alpha,\beta}$ with $o\left(  T_{n}\right)
=\omega^{\alpha+\beta\cdot n}.$ Hence $I_{b}\left(  \overset{_{\sim}}%
{T}_{\alpha,\beta},2^{n}\right)  \geq\omega^{\alpha+\beta\cdot n}.$ Thus
$I_{b}\left(  \overset{_{\sim}}{T}_{\alpha,\beta}\right)  =\sup_{K}%
I_{b}\left(  \overset{_{\sim}}{T}_{\alpha,\beta},K\right)  \geq\omega
^{\alpha+\beta\cdot\omega}.$ We conclude that $I_{b}\left(  \overset{_{\sim}%
}{T}_{\alpha,\beta}\right)  =\omega^{\alpha+\beta\cdot\omega}.$ As $I\left(
\overset{_{\sim}}{T}_{\alpha,\beta}\right)  \geq I_{b}\left(  \overset{_{\sim
}}{T}_{\alpha,\beta}\right)  \geq\omega^{\omega},$ it follows from
\cite[Corollary 5.13]{JO} that $I\left(  \overset{_{\sim}}{T}_{\alpha,\beta
}\right)  =I_{b}\left(  \overset{_{\sim}}{T}_{\alpha,\beta}\right)  .$
\end{proof}

For the final corollary, recall that the Schreier space $X_{\alpha}%
,\,\alpha<\omega_{1},$ is the completion of $c_{00}$ with respect to the norm
$\left\|  \cdot\right\|  _{\mathcal{S}_{\alpha}}.$

\begin{corollary}
Suppose $\alpha$ is a countable ordinal whose Cantor normal form is
$\omega^{\alpha_{1}}\cdot m_{1}+...+\omega^{\alpha_{k}}\cdot m_{k}.$ If
$\alpha_{k}$ is not a limit ordinal, then there exists a Banach space $X$ such
that $I\left(  X\right)  =\omega^{\alpha}.$
\end{corollary}

\begin{proof}
If $\alpha_{k}=0,$ then $\alpha$ is a successor ordinal and $\iota\left(
X_{\alpha-1}\right)  =\omega^{\alpha}$ (\cite{AJO}). If $\alpha_{k}$ is a
successor ordinal, let $\gamma=\omega^{\alpha_{1}}\cdot m_{1}+...+\omega
^{\alpha_{k}}\cdot\left(  m_{k}-1\right)  $ and $\eta=\omega^{\alpha_{k}-1}.$
By Theorem \ref{T14}, $I\left(  \overset{_{\sim}}{T}_{\gamma,\eta}\right)
=\omega^{\gamma+\eta\cdot\omega}=\omega^{\alpha}.$
\end{proof}

\begin{remark}
\emph{The following analog of Proposition \ref{L2} for the space }$X_{\alpha
},$\emph{ }$\alpha<\omega_{1},$\emph{ holds obviously: If }$\left(
x_{i}\right)  _{i=1}^{p}$\emph{ is a normalized block basis of the unit vector
basis }$\left(  e_{k}\right)  $\emph{ of }$X_{\alpha},$\emph{ and }$k_{i}%
=\max$\emph{ supp\thinspace}$x_{i},$\emph{ }$1\leq i\leq p,$\emph{ then
}$\left\|  \sum_{i=1}^{p}a_{i}x_{i}\right\|  \leq\left\|  \sum_{i=1}^{p}%
a_{i}e_{k_{i}}\right\|  $\emph{ for all }$\left(  a_{i}\right)  \in c_{00}%
.$\emph{ Therefore the arguments in Proposition \ref{L8} and Theorem \ref{T14}
can be used to compute }$I_{b}\left(  X_{\alpha}\right)  $\emph{ (with respect
to the basis }$\left(  e_{k}\right)  $\emph{).}
\end{remark}

\end{document}